\documentclass{tac}
\usepackage{amsmath,amssymb,amsfonts,latexsym,stmaryrd,diagrams,mathrsfs}
\usepackage{cmll,bbm}

\usepackage{mymacros}

\title[Simplicial sets inside cubical sets]
{Simplicial sets inside cubical sets}
\author{Thomas Streicher and Jonathan Weinberger}
\copyrightyear{2021}

\eaddress{streicher@mathematik.tu-darmstadt.de\CR weinberger@mathematik.tu-darmstadt.de}

\date{\today}

\address{Department of Mathematics, TU Darmstadt, Schlossgartenstrasse 7, 64289 Darmstadt, Germany.}

\thanks{
We are grateful to Christian Sattler for many discussions, pointing out mistakes as well as problems in previous versions, and finally for making \cite{Sat18}, v2, available on the arXiv. Furthermore, we thank the reviewer for valuable comments and suggestions.
}


\begin{document}

\maketitle
\begin{abstract}
As observed recently by various people the topos $\sSet$ of simplicial sets appears as essential subtopos of a
topos $\cSet$ of cubical sets, namely presheaves over the category $\FL$ of finite lattices and monotone maps
between them. The latter is a variant of the cubical model of type theory due to Cohen \emph{et al.}\ for the purpose of providing a model for a variant of type theory which validates Voevodsky's Univalence Axiom and has computational meaning. 

Our contribution consists in constructing in $\cSet$ a fibrant univalent universe for those types that are sheaves. This makes it possible to consider $\sSet$ as a submodel of $\cSet$ for univalent Martin-L\"of type theory.

Furthermore, we adress the question whether the type-theoretic Cisinski model structure considered on $\cSet$ coincides with the test model structure, the latter of which models the homotopy theory of spaces. We do not provide an answer to this open problem, but instead give a reformulation in terms of the adjoint functors at hand.
\end{abstract}

\section{Introduction}

As observed in \cite{HS94} intensional Martin-L\"of type theory should have a natural interpretation in weak $\infty$-groupoids. In the first decade of this millennium it was observed that simplicial sets are a possible implementation
of this idea and around 2006 Voevodsky proved that the universe in question validates the so-called \emph{Univalence Axiom}
(UA) which roughly speaking states that isomorphic types are propositionally equal, see \cite{KL12} for a detailed proof.

But adding a constant inhabiting the type expressing UA gives rise to a type theory lacking computational meaning.
To overcome this problem Coquand \emph{et al.}\ \cite{CCHM18} have developed a so-called \emph{Cubical Type Theory}
based on explicit box filling operations from which UA can be derived. Cubes are finite powers of an interval
object $\Ibb$ which itself is not a proper type. In \cite{CCHM18} this type theory is interpreted
in the topos of covariant presheaves over the category of finitely presented free de Morgan algebras. It is clear
that (at least) the standard intensional Martin-L\"of type theory fragment together with the Univalence Axiom can be interpreted in the presheaf category $\cSet$ over the site $\Box$ which
is the full subcategory of $\Poset$ on finite powers of the 2 element lattice $\Ibb$. This site is op-equivalent
to the algebraic theory of distributive lattices as observed in \cite{Spi16}.

As observed independently in \cite{KV17} the topos $\sSet$ appears as subtopos of $\cSet$
and actually as an essential subtopos, cf.~\cite{Sat18}.

Starting from a universe of fibrant cubical sets within $\cSet$ we will construct a universe for fibrant simplicial sets
within $\cSet$. As it turns out, this universe is itself fibrant and univalent.

\section{Simplicial sets inside cubical sets}

We write $\Delta$ for the full subcategory of $\Poset$ on finite ordinals greater $0$
and we write $\Box$ for the full subcategory of $\Poset$ on finite powers of $\Ibb$, 
the 2 element lattice. Presheaves over $\Delta$ are called \emph{simplicial sets} and 
presheaves over $\Box$ are called \emph{cubical sets}. We write $\sSet$ and $\cSet$
for the toposes of simplicial and cubical sets, respectively.

Kapulkin and Voevodsky have observed in \cite{KV17} that one may obtain
$\sSet$ as a subtopos of $\cSet$ in the following way. The nerve
functor $\Nv : \Cat \to \sSet$ is known to be full and faithful and
so is its restriction $u : \Box \to \sSet$ to the full subcategory
$\Box$ of $\Cat$. This functor $u$ induces an adjunction 
$u_! \dashv u^* : \sSet \to \cSet$ where $u^*(X) = \sSet(u(-),X)$ and 
$u_!$ is the left Kan extension of $u$ along $\Yon_\Box : \Box\to\cSet$.
It follows from general topos theoretic results that $u_! \dashv u^*$
exhibits $\sSet$ as a subtopos of $\cSet$ induced by the Grothendieck
topology $\Jcal$ consisting of those sieves in $\Box$ which are sent 
by $u$ to jointly epic families in $\sSet$.

A more direct proof of a stronger result has been found by Sattler
in \cite{Sat18} and independently by the authors of this paper based
on the well known fact that splitting idempotents in $\Box$ gives
rise to the category of finite lattices and monotone maps between
them. E.g.\ by restricting to subobjects of objects in $\Box$ one
obtains an equivalent small full subcategory $\FL$ of $\Poset$.
Thus $\cSet$ is equivalent to $\widehat{\FL} = \Set^{\FL^\op}$ for which
reason we write $\cSet$ for $\widehat{\FL}$. 

The inclusion functor $i : \Delta \to \FL$ induces an essential geometric morphism
$i_! \dashv i^* \dashv i_*$ which, moreover, is injective, \emph{i.e.}\ $i_*$ and
thus also $i_!$ is full and faithful.
The inverse image part $i^*$ restricts presheaves over $\FL$ to
presheaves over $\Delta$ (by precomposition with $i^\op$). The direct image
part $i_*$ is given by $i_*(X) = \sSet(\Nv(-),X)$ since $\Nv$ restricted to
$\FL$ is given by $i^* \circ \Yon_\FL$.
The cocontinuous functor $i_!$ is the left Kan extension of $\Yon_\FL \circ i$
along $\Yon_\Delta$. It sends $X \in \sSet$ to the colimit of
$\Delta{\downarrow}X \stackrel{\partial_0}{\to} \Delta \stackrel{i}{\to}
\FL \stackrel{\Yon_\FL}{\longrightarrow} \cSet$.

Although not needed in full detail later on, we now give a rather explicit description of the Grothendieck topology corresponding to the injective geometric
morphism $i^* \dashv i_*$:  $S \subseteq \Yon_\FL(L)$ 
is a cover iff $i^*S = i^*\Yon_\FL(L)$, \emph{i.e.}\ $S$ contains all chains in $L$, \emph{i.e.}\
all monotone maps $c : [n] \to L$.\footnote{Thus, a sieve $S \subseteq \Yon_\FL(\Ibb^n)$
covers iff for every maximal chain $C \subseteq \Ibb^n$ there is an idempotent
$r \in S$ whose image is $C$. As mentioned earlier, recall that finite lattices are the idempotent completion of the cube category $\Box$.} Obviously, such an $S$ contains all monotone 
maps to $\Ibb^n$ whose image is contained in $C$. Thus the collection of all 
monotone maps to $\Ibb^n$ whose image is contained in a (maximal) chain in $\Ibb^n$
is the least covering sieve for $\Ibb^n$.

\section{Type-theoretic model structures on simplicial and cubical sets}

The representable object $\Ibb$ in $\cSet$ and $\sSet$ induces a Cisinski model structure
on these presheaf toposes which is generated by open box inclusions as described in the
following definition. For a general account cf.~\cite{Cis19}, Sec.~2.4., Thm.~2.4.19.

\definition[Type-theoretic model structure on $\sSet$ and $\cSet$]
Let $\mathcal{E}$ be $\cSet$ or $\sSet$. The \emph{interval} $\Ibb$ is given by $\Yon_{\FL}(\Ibb)$ and 
$\Yon_{\Delta}(\Ibb)$, respectively.
The \emph{type-theoretic model structure on $\mathcal E$} is defined by taking as cofibrations the monomorphisms and as fibrations those maps which are weakly right orthogonal to all \emph{open box inclusions} $(\{\varepsilon\} \times X) \cup (\Ibb \times Y) \inc \Ibb \times X$ where $Y \inc X$
and $\varepsilon \in \{0,1\}$, 
\enddefinition

In fact, this Cisinski model structure on simplicial sets coincides with the well known Kan model structure
since by \cite{GJ}, Ch.~I, Prop.~4.2, open box inclusions and horn inclusions generate the same class of 
anodyne extensions. 

Furthermore, it is obtained just by restricting the type-theoretic model structure on cubical sets defined above to simplicial sets.

In this section, we present a proof of this fact. All of the results have previously been established by Sattler in \cite{Sat18}, cf.~Prop.~3.3 and Sec.~3.3. We recall them here for clarification of context for the later sections, in particular for direct use in the universe construction of Sec.~\ref{sec:univ}.

\proposition\label{prop:i-lower-fib}
The inclusion $i_*:\sSet \to \cSet$ preserves and reflects fibrations.
\endproposition

\proof
From the preservation properties of the sheafification functor $i^*$ it 
follows immediately that open box inclusions in $\sSet$ are precisely the
sheafifications of open box inclusions in $\cSet$.

Since both in $\cSet$ and $\sSet$ fibrations are those maps which are weakly right 
orthogonal to all open box inclusions a map $f$ in $\sSet$ is a
fibration iff $i_*f$ is a fibration in $\cSet$.
\endproof

Writing $\F$ for the class of fibrations in $\cSet$ the class of
Kan fibrations in $\sSet$ is given by $\F \cap \sSet$ (considering
$\sSet$ as full subcategory of $\cSet$ via $i_*$).

Sattler has pointed out to us an elegant argument that $i_*$ preserves and reflects weak equivalences between fibrant objects.

\proposition\label{prop:weqpresrefl}
For fibrant objects $A,B \in \sSet$ a map $f : A \to B$ is a weak equivalence in $\sSet$ 
iff $i_*f$ is a weak equivalence in $\cSet$.
\endproposition

\proof
By \cite{Cis19}, Prop.~2.4.26, since both $\sSet$ and $\cSet$ are Cisinski model categories weak equivalences between fibrant objects are just homotopy equivalences 
and these are preserved by $i^*$ and $i_*$ since these functors preserve $\Ibb$ and finite products.
\endproof

\theorem
The adjunction $i^* \dashv i_*$ is a Quillen adjunction. Moreover, $i^*$ preserves weak equivalences between arbitrary objects.
\endtheorem

\proof
By Proposition \ref{prop:i-lower-fib} the functor $i_*$ preserves fibrations. As $i^*$ is in turn a right adjoint it preserves monomorphisms. Thus $i^* \dashv i_*$ is a Quillen adjunction.

As every cubical set is cofibrant, by Ken Brown's Lemma $i^*$ preserves weak equivalences.
\endproof

Even more is true, namely $i_! \dashv i^*$ is a Quillen adjunction as well, as shown by Sattler in \cite{Sat18}, Sec.~3.3. Hence, in particular $i^*$ also preserves fibrations.

\theorem\label{firstqadj}
The adjunction $i_! \dashv i^*$ is a Quillen adjunction.
\endtheorem

Alas, it is not known whether $i^* \dashv i_*$ is a Quillen \emph{equivalence}. For that purpose one would have to show that for
fibrant $B \in \sSet$ and arbitrary $A \in \cSet$ a map $f : i^*A \to B$ is a
weak equivalence in $\sSet$ iff the transpose $\check{f} : A \to i_*B$ is a weak equivalence
in $\cSet$.

\section{Universes in cubical sets}\label{sec:univ}

In general, given a finite limit category $\mathbf{C}$ small full subfibrations of the fundamental (codomain) fibration $P_\mathbf{C}: \mathbf{C}^{\mathbf 2} \to \mathbf{C}$ are given by pullback-stable classes $\mathcal{S}$ of morphisms in $\mathbf{C}$ admitting a \emph{generic family}, \emph{i.e.}~a map $\pi \in \mathcal S$ s.t.~every map $f \in \mathcal{S}$ arises as pullback of $\pi$, cf.~\cite{Str20}.

Given a Grothendieck universe $\Ucal$ this induces a universe
\`a la Yoneda $\pi:E \to U$ in $\widehat{\FL}$  which is generic 
for the class of $\Ucal$-small maps in $\widehat{\FL}$ \cite{Str05}.

Now, as described in \cite[Sect.~9]{GS17} there is a universe $\pi_c : E_c \to U_c$
generic\footnote{\emph{i.e.} all $\Ucal$-small fibrations can be obtained as pullback of the
generic one in a typically non-unique way} for $\Ucal$-small fibrations in $\cSet$ such that $U_c$
is fibrant. Moreover, the universe $\pi_c$ indeed arises as a subuniverse of $\pi$ via a map $u_c: U_c \to U$. Note that as described in \cite{GS17} $U_c(L)$ does not simply consist of 
$\Ucal$-small fibrations over $\Yon_\FL(L)$ but rather such fibrations together with a
functorial choice of fillers which are forgotten by $u_c$.

Recall from \cite{KL12}, Def.~3.2.10, the definition of a univalent universe. By Kripke--Joyal translation, this amounts to the following. A universe $\pi: E \to U$ is \emph{univalent} if and only if for any $f,g: I \to U$ from $f^*\pi  \sim g^*\pi$ it follows that $f \sim g$, \emph{i.e.}~classifying maps are unique up to homotopy.

We now give the construction of a fibrant univalent universe in cubical sets which is 
generic for small fibrations which are families of sheaves. Recall from \cite{Str05} 
that a map $a : A \to I$ is a family of sheaves iff the naturality square
\begin{diagram}[small]
A \SEpbk & \rTo^{\eta_A} & i_*i^*A \\
\dTo^{a} &               & \dTo_{i_*i^*a} \\ 
I        & \rTo_{\eta_I} & i_*i^*I 
\end{diagram}
is a pullback. As follows from \cite{Str05} $i_*i^*\pi_c$ is a generic family for maps in $\cSet$ which are small fibrations and  families of sheaves  but, alas, this universe is not univalent. In
the rest of this section we construct a univalent universe from this
making, however, use of the univalent universe $\pi_s : E_s \to U_s$ 
in $\sSet$ as constructed in \cite{KL12}.

\theorem\label{thm:univ}
A univalent small fibration generic for small fibrations which are also families 
of sheaves can be obtained by pulling back $i_*i^*\pi_c$ along the homotopy
equalizer of $i_*e \circ i_*p$ and $\id_{i_*U_s}$ where $e$ and $p$ are maps
such that both squares in
\begin{diagram}
	E_s \SEpbk & \rTo & i^*E_c & &  i^*E_c \SEpbk & \rTo & E_s \\
	\dTo^{\pi_s} & & \dTo_{i^*\pi_c} &  &\dTo^{i^*\pi_c} & & \dTo_{\pi_s} \\  
	U_s & \rTo_e & i^*U_c & & i^*U_c & \rTo_p & U_s
\end{diagram}
are pullbacks.
\endtheorem

\proof
We start with two universal fibrations, namely $\pi_s: E_s \to U_s$ in $\sSet$ classifying for small Kan fibrations, and $\pi_c: E_c \to U_c$ in $\cSet$ classifying for small fibrations in the type-theoretic model structure.

By Thm.~\ref{firstqadj} the functor $i^*$ preserves fibrations for which reason $i^*\pi_c$ is a fibration. Thus, there is a map $p:i^* U_c \to U_s$ fitting into a pullback square:
\begin{diagram}
 i^*E_c \SEpbk & \rTo & E_s \\
 \dTo^{i^*\pi_c} & & \dTo_{\pi_s} \\  
i^*U_c & \rTo_p & U_s
\end{diagram}
On the other hand, small fibrations which are families of sheaves are in the essential image of $i^*$ when restricted to fibrations in $\cSet$. Thus, we obtain that in turn $\pi_s$ arises as a pullback of $i^*\pi_c$, \emph{i.e.}~there is a map $e:U_s \to i^*U_c$ such that we have a pullback square:
\begin{diagram}
	E_s \SEpbk & \rTo & i^*E_c \\
	\dTo^{\pi_s} & & \dTo_{i^*\pi_c}  \\  
	U_s & \rTo_e & i^*U_c 
\end{diagram}
Since the universe $\pi_s$ is univalent, it is \emph{classifying} up to homotopy, rather than merely generic. Thus, pasting the two pullback squares together we obtain that $U_s$ classifies itself through the composite $p \circ e: U_s \to U_s$. By uniqueness up to homotopy, $p \circ e$ must then be homotopic to $\id_{U_s}$:
\begin{diagram}
	E_s \SEpbk & \rTo & i^*E_c \SEpbk & \rTo & E_s \\
	\dTo^{\pi_s} & & \dTo_{i^*\pi_c} & & \dTo_{\pi_s} \\  
	U_s & \rTo_e & i^*U_c & \rTo_p & U_s
\end{diagram}

 Since $i_*$ preserves fibrations by Prop.~\ref{prop:i-lower-fib}, pullbacks and $\sim$ (homotopy) for maps between fibrant objects we have
\begin{diagram}
i_*E_s \SEpbk & \rTo & i_*i^*E_c \SEpbk & \rTo & i_*E_s \\
\dTo^{i_*\pi_s} & & \dTo_{i_*i^*\pi_c} & & \dTo_{i_*\pi_s} \\  
i_*U_s & \rTo_{i_*e} & i_*i^*U_c & \rTo_{i_*p} & i_*U_s
\end{diagram}
with $i_*p \circ i_*e = i_*(p \circ e) \sim i_*(\id_{U_s}) = \id_{i_*U_s}$. We want to argue that $i_* \pi_s$ is a univalent universe for small fibrations that are families of sheaves. 

For this purpose, consider a small fibration $p:A \to I$ together with a diagram:
\begin{diagram}
	A\SEpbk &  \pile{\rTo^a \\ \rTo_b} & i_*E_s \\
	\dTo^{p} & (1)  & \dTo_{i_* \pi_s}  \\  
	I &  \pile{\rTo^f\\ \rTo_g} & i_*U_s
\end{diagram}
This factors as follows:
\begin{diagram}
	A & \rTo^{\eta_A} & i_*i^*A \SEpbk &  \pile{\rTo^{a'} \\ \rTo_{b'}} & i_* E_s \\
	\dTo^{p} & \hspace{-1cm} (2) & \dTo^{i_*i^*p} & (3) & \dTo_{i_* \pi_s} \\  
	I &   \rTo_{\eta_I} & i_*i^*I &  \pile{\rTo^{f'}\\ \rTo_{g'}}  & i_* U_s
\end{diagram}
Square (3) is the image of Square (1) under the reflection $i_*i^*$, so it is indeed a pullback. By the Pullback Lemma, then Square (2) is a pullback, too: 
\begin{diagram}
	A \SEpbk & \rTo^{\eta_A} &  i_*i^*A \\
	\dTo^{p} & & \dTo_{i_*i^*\pi_s}  \\  
	I &   \rTo_{\eta_I} &  i_*i^* I 
\end{diagram}
This means, $p:A \to I$ is a family of sheaves, proving one part of the claim.
Now, the functor $i^*$ preserves pullbacks, thus it maps Square (3) to
\begin{diagram}
	i^*A \SEpbk &  \pile{\rTo^{i^*a'} \\ \rTo_{i^*b'}} & E_s \\
	\dTo^{i^*p} &   & \dTo_{\pi_s}  \\  
	i^*I &  \pile{\rTo^{i^*f'}\\ \rTo_{i^*g'}} & U_s
\end{diagram}
where we assume for sake of simplicity that $i^*i_*$ is the identity functor.

Univalence of $\pi_s$ implies $i^*f' \sim i^*g'$. Since $\sSet$ arises as a full subcategory of $\cSet$, there exist maps $\widetilde{f}$ and $\widetilde{g}$ in $\sSet$ such that $f' = i_* \widetilde{f}$ and $g'  = i_* \widetilde{g}$, hence $\widetilde{f} \sim \widetilde{g}$. Since $i_*$ preserves homotopy between maps, $f' = i_* \widetilde{f} \sim i_* \widetilde{g} = g'$. Finally, this implies $f \sim g$. 

Thus, $i_* \pi_s$ is a fibrant univalent universe for small fibrations that are families of sheaves.

Now, considering the diagram
\begin{diagram}
i_* i^*E_c  \SEpbk & \rTo & 	i_* E_s \SEpbk & \rTo &i_* i^*E_c \\
	\dTo^{i_* i^*\pi_c} & & \dTo^{i_* \pi_s} & & 	\dTo^{i_* i^*\pi_c} \\  
i_* i^*U_c  & \rTo_{i_*p} & i_*U_s & \rTo_{i_*e} & i_* i^*U_c 
\end{diagram}
we find that $i_*i^*\pi_c$ is generic for small fibrations which are families of
sheaves since $i_*\pi_s$ is a classifying fibration and the right-hand square commutes.

As argued before, the maps $i_*e$ and $i_*p$ form a homotopy section-retraction pair. Hence, $i_*e$ is a homotopy equalizer of $i_*e \circ i_*p$ and $\id_{U_s}$. Pulling back $i_* i^* \pi_c$ along this homotopy equalizer yields a univalent subuniverse of the constructive universe $i_* i^* \pi_c$ in a constructive way. But as shown above $i^*i_*\pi_c$ is nonconstructively equivalent to $i_*\pi_c$ which is a universe in a nonconstructive way.
\endproof

\section{Does the type-theoretic model structure on $\cSet$ coincide with the test model structure?}

Having discussed how the type-theoretic model structure coincides with the standard Kan model structure on simplicial sets the following question arises: Does the type-theoretic model structure on \emph{cubical sets} also present the (standard) homotopy theory of spaces? To this date, the question remains unanswered. In the section at hand, we give a formulation of this problem in terms of \emph{test model structures}, using the Quillen adjunctions discussed previously.

Introduced by Grothendieck \cite{pstacks}, \emph{test categories} admit a certain model structure on their presheaf category whose homotopy category is equivalent to the standard homotopy category of spaces. A general comprehensive theory of test model structures has been developed by Cisinski in his thesis \cite{Cis06}. Further accounts are given in Jardine \cite{Jar06} and Maltsiniotis \cite{Mal05}.

In particular the simplex category $\Delta$ is a test category (\cite{Mal05}, Prop.~1.5.13), as are most of the familiar cube categories, cf. \cite{Cis06}, Ch.~8, \cite{Jar06}, Sec.~8, and \cite{BM17}.

For a concise recollection of the notions of test category and test model structure cf.~\cite{Jar06}, Sec.~2.

In our setting, we can ask if the type-theoretic model structure \emph{on cubical sets} coincides with the test model structure. Since this is still an open problem, we are not providing an answer, but rather an interesting reformulation of the problem. Namely, the type-theoretic model structure on $\cSet$ coincides with the test model structure if and only if all the components of the counit of $i_! \dashv i^*$ are weak equivalences.

We begin our discussion by noting that the inclusion functor $i : \Delta \inc \FL$ is aspherical in the sense of
\cite{Mal05}, Sec.~Def.~1.1.2, \emph{i.e.}~$\Nv(i \downarrow L)$ is contractible in the Kan model structure on $\sSet$, for all $L \in \FL$. This follows since every comma category $i \downarrow L$ is connected. 
Thus by \cite[Th.~1.2.9]{Mal05} the functor $i^* : \cSet \to \sSet$ preserves 
and reflects weak equivalences of the respective test model structures. Since 
$i^*$ also preserves monos the adjunction $i^* \dashv i_*$ is a Quillen 
equivalence between $\cSet$ and $\sSet$ endowed the respective test model
structures.

Let $\varepsilon_X : i_!i^*X \to X$ be the counit of $i_!\dashv i^*$ and
$\eta_X : X \to i_*i^*X$ be the unit of $i^* \dashv i_*$. Both maps are sent to isos by $i^*$ and thus are weak
equivalences w.r.t.\ the test model structure on $\cSet$.

We know that both $i_! \dashv i^*$ and $i^* \dashv i_*$ are Quillen pairs 
when $\cSet$ is endowed with the type-theoretic model structure.
Thus, if the type-theoretic model structure on $\cSet$ coincides with 
the test model structure then all $\eta_X : X \to i_*i^*X$ and 
$\varepsilon_X : i_!i^*X \to X$ are weak equivalences w.r.t.\ the 
type-theoretic model structure on $\cSet$.
But if all $\varepsilon_X$ are weak equivalences w.r.t.\ the type-theoretic model structure then it coincides with the test model structure which can
be seen as follows. Suppose $m : Y \to X$ in $\cSet$ is an anodyne cofibration w.r.t.\ the 
test model structure then $i^*m$ is an anodyne cofibration in $\sSet$ from 
which it follows that $i_!i^*m$ is an anodyne cofibration w.r.t.\ the type-theoretic model structure on $\cSet$. But since
\begin{diagram}[small]
i_!i^*Y & \rTo^{\varepsilon_Y} & Y \\
\dTo^{i_!i^*m} & & \dTo_m \\
i_!i^*X & \rTo_{\varepsilon_X} & X \\
\end{diagram}
commutes it follows by the 2-out-of-3 property for weak equivalences
that $m$ is a weak equivalence and thus an anodyne cofibration w.r.t.\ 
the type-theoretic model structure on $\cSet$.

Thus, summarizing the above considerations we conclude that the 
type-theoretic and the test model structure on $\cSet$ coincide
if and only if all $\varepsilon_X : i_!i^*X \to X$ are weak equivalences
in the type-theoretic model structure on $\cSet$. Alas, we do not know 
whether this is the case in general.

\section{Conclusion}

Using the fact that $\sSet$ is an essential subtopos of $\cSet$ we have constructed a fibrant univalent universe inside $\cSet$ which is generic for small families of sheaves, \emph{i.e.}~simplicial sets.

However, this construction makes use of the univalent universe inside $\sSet$. Formally speaking this universe can be constructed in the internal language of $\cSet$ but only at the price of importing the inconstructive universe $\pi_s$ from $\sSet$ via $i_*$.

Nevertheless, this may still model an extension of the cubical type theory of \cite{CCHM18} providing a univalent universe for small simplicial sets the precise formulation of which we leave for future work.

\nocite{*}
\bibliographystyle{amsalpha}
\bibliography{cusi}{}

\providecommand{\bysame}{\leavevmode\hbox to3em{\hrulefill}\thinspace}
\providecommand{\MR}{\relax\ifhmode\unskip\space\fi MR }
\providecommand{\MRhref}[2]{%
  \href{http://www.ams.org/mathscinet-getitem?mr=#1}{#2}
}
\providecommand{\href}[2]{#2}
\begin{thebibliography}{CCHM18}

\bibitem[BM17]{BM17}
U.~Buchholtz and E.~Morehouse, \emph{Varieties of cubical sets}, Relational and
  Algebraic Methods in Computer Science (Cham) (Peter H{\"o}fner, Damien Pous,
  and Georg Struth, eds.), Springer International Publishing, 2017,
  \href{https://link.springer.com/chapter/10.1007/978-3-319-57418-9_5}{doi:10.1007/978-3-319-57418-9{\_}5},
  pp.~77--92.

\bibitem[CCHM18]{CCHM18}
C.~Cohen, T.~Coquand, S.~Huber, and A.~M{\"o}rtberg, \emph{{C}ubical {T}ype
  {T}heory: A constructive interpretation of the {U}nivalence {A}xiom}, 21st
  International Conference on Types for Proofs and Programs (TYPES 2015)
  (Dagstuhl, Germany) (Tarmo Uustalu, ed.), Leibniz International Proceedings
  in Informatics (LIPIcs), vol.~69, Schloss Dagstuhl--Leibniz-Zentrum fuer
  Informatik, 2018,
  \href{https://drops.dagstuhl.de/opus/volltexte/2018/8475/}{doi:10.4230/LIPIcs.TYPES.2015.5},
  pp.~5:1--5:34.

\bibitem[Cis06]{Cis06}
D.-C. Cisinski, \emph{Les pr\'efaisceaux comme mod\`eles des types
  d'homotopie}, Ast\'erisque, no. 308, Soci\'et\'e math\'ematique de France,
  2006 (fr),
  \href{http://www.numdam.org/item/AST_2006__308__R1_0/}{doi:10.24033/ast.715}.
  \MR{2294028}

\bibitem[Cis19]{Cis19}
D.-C. Cisinski, \emph{Higher categories and homotopical algebra}, Cambridge
  University Press, 2019,
  \href{https://doi.org/10.1017/9781108588737}{doi:10.1017/9781108588737}.

\bibitem[GJ09]{GJ}
P.~G. Goerss and J.~F. Jardine, \emph{Simplicial homotopy theory}, Springer,
  2009,
  \href{https://link.springer.com/book/10.1007/978-3-0346-0189-4}{doi:10.1007/978-3-0346-0189-4}.

\bibitem[Gro83]{pstacks}
A.~Grothendieck, \emph{Pursuing stacks (original: {\`A} la poursuite des
  champs)}, 1983,
  \url{https://thescrivener.github.io/PursuingStacks/ps-online.pdf}.

\bibitem[GS17]{GS17}
N.~Gambino and C.~Sattler, \emph{The {F}robenius condition, right properness,
  and uniform fibrations}, Journal of Pure and Applied Algebra \textbf{221}
  (2017), no.~12, 3027 -- 3068,
  \href{https://doi.org/10.1016/j.jpaa.2017.02.013}{doi:10.1016/j.jpaa.2017.02.013}.

\bibitem[GZ67]{GZ67}
P.~Gabriel and G.~Zisman, \emph{Calculus of fractions and homotopy theory},
  Springer, 1967,
  \href{https://doi.org/10.1112/jlms/s1-44.1.382a}{doi:10.1007/978-3-642-85844-4}.

\bibitem[HS94]{HS94}
M.~{Hofmann} and T.~{Streicher}, \emph{The groupoid model refutes uniqueness of
  identity proofs}, Proceedings Ninth Annual IEEE Symposium on Logic in
  Computer Science, July 1994,
  \href{https://doi.ieeecomputersociety.org/10.1109/LICS.1994.316071}{doi:10.1109/LICS.1994.316071},
  pp.~208--212.

\bibitem[Jar06]{Jar06}
J.~F. Jardine, \emph{Categorical homotopy theory}, Homology Homotopy Appl.
  \textbf{8} (2006), no.~1, 71--144,
  \href{https://www.intlpress.com/site/pub/pages/journals/items/hha/content/vols/0008/0001/a003/index.php}{doi:10.4310/HHA.2006.v8.n1.a3}.

\bibitem[KL12]{KL12}
K.~Kapulkin and P.~Lumsdaine, \emph{The simplicial model of univalent
  foundations (after {V}oevodsky)}, Journal of the European Mathematical
  Society (2012), forthcoming,
  \href{https://arxiv.org/abs/1211.2851}{arXiv:1211.2851}.

\bibitem[KV20]{KV17}
K.~Kapulkin and V.~Voevodsky, \emph{A cubical approach to straightening},
  Journal of Topology \textbf{13} (2020), no.~4, 1682--1700,
  \href{https://doi.org/10.1112/topo.12173}{doi:10.1112/topo.12173}.

\bibitem[Mal05]{Mal05}
G.~Maltsiniotis, \emph{La th{\'e}orie de l'homotopie de {G}rothendieck},
  Soci{\'e}t{\'e} Math{\'e}matique de France, 2005,
  \href{https://smf.emath.fr/publications/la-theorie-de-lhomotopie-de-grothendieck}{doi:10.24033/ast.689},
  updated and expanded version:
  \url{https://webusers.imj-prg.fr/~georges.maltsiniotis/ps/prstnew.pdf}.

\bibitem[Sat18]{Sat18}
C.~Sattler, \emph{Idempotent completion of cubes in posets},
  \href{https://arxiv.org/abs/1805.04126v2}{arXiv:1805.04126v2}, 2018.

\bibitem[Spi16]{Spi16}
B.~Spitters, \emph{Cubical sets and the topological topos},
  \href{https://arxiv.org/abs/1610.05270}{arXiv:1610.05270}, 2016.

\bibitem[{Str}05]{Str05}
T.~{Streicher}, \emph{{U}niverses in {T}oposes}, From Sets and Types to
  Topology and Analysis, Towards Practicable Foundations for Constructive
  Mathematics \textbf{48} (2005), 78--90,
  \href{https://oxford.universitypressscholarship.com/view/10.1093/acprof:oso/9780198566519.001.0001/acprof-9780198566519-chapter-5}{doi:10.1093/acprof:oso/9780198566519.003.0005
  }.

\bibitem[Str20]{Str20}
T.~Streicher, \emph{Fibred categories \`{a} la {J}ean {B}\'{e}nabou},
  \href{https://arxiv.org/abs/1801.02927}{arXiv:1801.02927}, 2020.

\end{thebibliography}

\end{document}